\documentclass[runningheads]{llncs}
\usepackage{graphicx}
\usepackage{hyperref}
\hypersetup{
    colorlinks=true,
    linkcolor=blue,
    filecolor=magenta,      
    urlcolor=blue,}
\usepackage{amsfonts, amsfonts, amssymb, amsmath}
\usepackage{color}
\newcommand{\g}{\mathfrak{g}}
\newcommand{\R}{\mathbb{R}}
\newcommand{\Ad}{\textrm{Ad}}
\newcommand{\ad}{\textrm{ad}}
\newcommand{\se}{\mathfrak{se}}

\begin{document}
\title{A reduced parallel transport equation on Lie Groups with a left-invariant metric}
\titlerunning{Parallel Transport on Lie Groups}
%
\author{Nicolas Guigui\orcidID{0000-0002-7901-0732} \and
Xavier Pennec}
\authorrunning{N. Guigui and X. Pennec}
%
\institute{Université Côte d'Azur, Inria Epione project team, France}
\maketitle              
\begin{abstract}
This paper presents a derivation of the parallel transport equation expressed in the Lie algebra of a Lie group endowed with a left-invariant metric.
The use of this equation is exemplified on the group of rigid body motions $SE(3)$, using basic numerical integration schemes, and compared to the pole ladder algorithm. This results in a stable and efficient implementation of parallel transport. The implementation leverages the python package \url{geomstats} and is available online.
\end{abstract}
\section{Introduction}
Lie groups are ubiquitous in geometry, physics and many application domains such as robotics \cite{barrau_invariant_2017}, medical imaging \cite{lorenzi_efficient_2014} or computer vision \cite{hauberg_unscented_2013}, giving rise to a prolific research avenue.
Structure preserving numerical methods have demonstrated significant qualitative and quantitative improvements over extrinsic methods \cite{iserles_lie-group_2005}. 
Moreover, machine learning \cite{barbaresco_lie_2020} and optimisation methods \cite{journee_gradient-optimization_2007} are being developed to deal with Lie group data.

In this context, parallel transport is a natural tool to define statistical models and optimisation procedures, such as the geodesic or spline regression \cite{kim_smoothing_2020,nava-yazdani_geodesic_2020}, or to normalise data represented by tangent vectors \cite{yair_parallel_2019,brooks_riemannian_2019}.

Different geometric structures are compatible with the group structure, such as its canonical Cartan connection, whose geodesics are one-parameter subgroups, or left-invariant Riemannian metrics.
In this work we focus on the latter case, that is fundamental in geometric mechanics \cite{kolev_lie_2004} and has been studied in depth since the foundational papers of Arnold \cite{arnold_sur_1966} and Milnor \cite{milnor_curvatures_1976}. The fundamental idea of Euler-Poincarré reduction is that the geodesic equation can be expressed entirely in the Lie algebra thanks to the symmetry of left-invariance~\cite{marsden_mechanical_2009}, alleviating the burden of coordinate charts.

However, to the best of our knowledge, there is no literature on a similar treatment of the parallel transport equation. We present here a derivation of the parallel transport equation expressed in the Lie algebra of a Lie group endowed with a left-invariant metric.
We exemplify the use of this equation on the group of rigid body motions $SE(3)$, using common numerical integration schemes, and compare it to the pole ladder approximation algorithm. This results in a stable and efficient implementation of parallel transport. The implementation leverages the python package \url{geomstats} and is available online at \url{http://geomstats.ai}.

In section~\ref{sec:notation}, we give the general notations and recall some basic facts from Lie group theory. Then  we derive algebraic expressions of the Levi-Civita connection associated to the left-invariant metric in section~\ref{sec:metric}. The equation of parallel transport is deduced from this expression  and its integration is exemplified in section~\ref{sec:equation}.

\section{Notations}
\label{sec:notation}
Let $G$ be a lie group of (finite) dimension $n$. Let $e$ be its identity element, $\g = T_eG$ be its tangent space at $e$, and for any $g \in G$, let $L_g: h\in G \mapsto gh$ denote the left-translation map, and $dL_g$ its differential map. Let $\g^L$ be the Lie algebra of left-invariant vector fields of $G$: $X \in \g^L \iff \forall g \in G, X|_g = dL_g X_e$.

$\g$ and $\g^L$ are in one-to-one correspondence, and we will write $\tilde x$ the left-invariant field generated by $x \in \g$: $\forall g \in G$, $\tilde x_g = dL_g x$.
The bracket defined on $\g$ by $[x, y] = [\Tilde{x}, \Tilde{y}]_e$ turns $\g$ into a Lie algebra that is isomorphic to $g^L$. One can also check that this bracket coincides with the adjoint map defined by $\ad_x(y) = d_e(g\mapsto \Ad_g y)$, where $Ad_g=d_e(h \mapsto ghg^{-1})$. For a matrix group, it is the commutator.

Let $(e_1, \ldots, e_n)$ be an orthonormal basis of $\g$, and the associated left-invariant vector fields $X_i^L = \Tilde{e_i} = g \mapsto dL_g e_i$. As $dL_g$ is an isomorphism, $(X_1^L|_g, \ldots, X_n^L|_g)$ form a basis of $T_gG$ for any $g \in G$, so one can write $X|_g = f^i(g)X_i^L|_g$ where for $i=1,\ldots,n$, $g\mapsto f^i(g)$ is a smooth real-valued function on $G$. Any vector field on $G$ can thus be expressed as a linear combination of the $X_i^L$ with function coefficients.

Finally, let $\theta$ be the Maurer-Cartan form defined on $G$ by:
\begin{equation}
    \forall g \in G, \forall v \in T_gG, \; \theta|_g(v) = (dL_g)^{-1} v \in \g
\end{equation}
It is a $\g$-valued 1-form and for a vector field $X$ on $G$ we write $\theta(X)|_g = \theta|_g(X|_g)$ to simplify the notations.

\section{Left-invariant metric and connection}
\label{sec:metric}
A Riemannian metric $\langle\cdot,\cdot\rangle$ on $G$ is called left-invariant if the differential map of the left translation is an isometry between tangent spaces, that is
\begin{equation*}
   \forall g,h \in G, \forall u,v \in T_gG, \;\; \langle u,v\rangle_g = \langle dL_h u, dL_h v\rangle_{hg}.
\end{equation*}
It is thus uniquely determined by an inner product on the tangent space at the identity $T_eG =\g$ of $G$. 
Furthermore, the metric dual to the adjoint map is defined such that
\begin{equation}
    \forall a,b,c \in \g, \langle\ad_a^*(b), c\rangle = \langle b, \ad_a(c)\rangle = \langle [a, c], b\rangle.
\end{equation}
As the bracket can be computed explicitly in the Lie algebra, so can $\ad^*$ thanks to the orthonormal basis of $\g$. Now let $\nabla$ be the Levi-Civita connection associated to the metric. It is also left-invariant and can be characterised by a bi-linear form on $\g$ that verifies \cite{pennec_exponential_2012,gallier_differential_2020}:
\begin{equation}
    \forall x,y \in \g, \;\; \alpha(x, y) := (\nabla_{\Tilde{x}}\Tilde{y})_e = \frac{1}{2}\big([x, y] - \ad_x^*(y) - \ad_y^*(x)\big)
    \label{eq:alpha}
\end{equation}
Indeed by the left-invariance, for two left-invariant vector fields $X=\tilde x, Y=\tilde y~\in~\g^L$, the map $g\mapsto \langle X, Y\rangle_g$ is constant, so for any vector field $Z=\tilde z$ we have $Z(\langle X,Y\rangle)=0$. Kozsul formula thus becomes
\begin{align}
        2\langle\nabla_X Y, Z\rangle &= \langle [X, Y], Z\rangle - \langle [Y, Z], X\rangle - \langle [X, Z], Y\rangle\label{eq:kozsul}\\
        2\langle\nabla_X Y, Z\rangle_e    &= \langle [x, y], z\rangle_e - \langle\ad_y(z), x\rangle_e - \langle\ad_x(z), y\rangle_e\nonumber\\
        2 \langle\alpha(x, y), z\rangle_e&= \langle [x, y], z\rangle_e- \langle\ad_y^*(x), z\rangle_e - \langle\ad_x^*(y), z\rangle_e.\nonumber
\end{align}
Note however that this formula is only valid for left-invariant vector fields. We will now generalise to any vector fields defined along a smooth curve on $G$, using the left-invariant basis ($X_1^L, \ldots, X_n^L$).

Let $\gamma:[0,1]\rightarrow G$ be a smooth curve, and $Y$ a vector field defined along $\gamma$.
Write $Y = g^i X_i^L$, $\dot{\gamma} = f^i X_i^L$. Let's also define the \textit{left-angular velocities} $\omega(t)= \theta|_{\gamma(t)} \dot{\gamma}(t) = (f^i\circ \gamma)(t) e_i \in \g$ and $\zeta(t) = \theta(Y)|_{\gamma(t)} = (g^j \circ \gamma)(t) e_j \in \g$. Then the covariant derivative of $Y$ along $\gamma$ is
\begin{align*}
    \nabla_{\dot \gamma(t)}Y &= (f^i \circ \gamma)(t) \nabla_{X_{i}^{L}}\big( g^i X_i^L \big) \\
        &= (f^i \circ \gamma)(t) X_i^L(g^j) X_j^L + (f^i \circ \gamma)(t) (g^j\circ \gamma)(t) (\nabla_{X_{i}^{L}}{X_j^L})_{\gamma(t)}\\
    dL_{\gamma(t)}^{-1} \nabla_{\dot \gamma(t)}Y &= (f^i \circ \gamma)(t) X_i^L(g^j) e_j + (f^i \circ \gamma)(t) (g^j\circ \gamma)(t) dL_{\gamma(t)}^{-1}(\nabla_{X_{i}^{L}}{X_j^L})_{\gamma(t)} \\
                &= (f^i \circ \gamma)(t) X_i^L(g^j) e_j + (f^i \circ \gamma)(t) (g^j \circ \gamma)(t)\nabla_{e_i}e_j
\end{align*}
where Leibniz formula and the invariance of the connection is used in $(\nabla_{X_{i}^{L}}{X_j^L}) = dL_{\gamma(t)} \nabla_{e_i}e_j$. Therefore for $k=1..n$
\begin{align}
   \langle dL_{\gamma(t)}^{-1} \nabla_{\dot \gamma(t)}Y, e_k\rangle &= (f^i \circ \gamma)(t) X_i^L(g^j) \langle e_j, e_k\rangle \nonumber \\
   &\quad + (f^i \circ \gamma)(t) (g^j\circ \gamma)(t) \langle\nabla_{e_i}e_j,e_k\rangle
\end{align}
but on one hand
\begin{align}
    \zeta(t) &= \theta(Y)|_{\gamma(t)} = \theta|_{\gamma(t)}\big(((g^j\circ \gamma)(t) X_j^L|_{\gamma(t)})\big) \nonumber\\
         &= (g^j\circ \gamma)(t) e_j \\
    \dot{\zeta}(t) & = (g^j\circ \gamma)'(t) e_j = d_{\gamma(t)}g^j \dot{\gamma}(t) e_j \nonumber\\
               &= d_{\gamma(t)} g^j \Big( (f^i \circ \gamma)(t) X_i^L|_{\gamma(t)}\Big) e_j \nonumber\\
               &= (f^i \circ \gamma)(t) d_{\gamma(t)}g^j X_i^L|_{\gamma(t)} e_j \nonumber\\
               &= (f^i \circ \gamma)(t) X_i^L(g^j) e_j
\end{align}
and on the other hand, using \eqref{eq:kozsul}:
\begin{align}
   (f^i \circ \gamma)(g^j\circ \gamma) \langle\nabla_{e_i}e_j,e_k\rangle  &= \frac{1}{2} (f^i \circ \gamma)(g^j\circ \gamma)(\langle [e_i, e_j], e_k\rangle \nonumber\\
   &\qquad - \langle [e_j, e_k], e_i\rangle - \langle [e_i, e_k], e_j\rangle) \nonumber\\
     &= \frac{1}{2} ( \langle [(f^i \circ \gamma) e_i, (g^j \circ \gamma) e_j], e_k\rangle \nonumber\\
     &\qquad - \langle [(g^j \circ \gamma) e_j, e_k], (f^i \circ \gamma) e_i\rangle \nonumber\\
     &\qquad - \langle [(f^i \circ \gamma) e_i, e_k], (g^j \circ \gamma) e_j\rangle) \nonumber\\
     &= \frac{1}{2}([\omega, \zeta] - \ad_{\omega}^*\zeta - \ad_{\zeta}^*\omega) = \alpha(\omega, \zeta)
\end{align}
Thus, we obtain an algebraic expression for the covariant derivative of any vector field $Y$ along a smooth curve $\gamma$. It will be the main ingredient of this paper.
\begin{equation}
\label{eq:connection}
    dL_{\gamma(t)}^{-1}\nabla_{\dot{\gamma}(t)}Y(t) = \dot{\zeta}(t) + \alpha(\omega(t), \zeta(t))
\end{equation}
A similar expression can be found in \cite{arnold_sur_1966,gay-balmaz_invariant_2012}. 
As all the variables of the right-hand side are defined in $\g$, they can be computed with matrix operations and an orthonormal basis.

\section{Parallel Transport}
\label{sec:equation}

We now focus on two particular cases of $\eqref{eq:connection}$ to derive the equations of geodesics and of parallel transport along a curve.

\subsection{Geodesic equation}
The first particular case is for $Y(t)=\dot \gamma(t)$. It is then straightforward to deduce from \eqref{eq:connection} the Euler-Poincarré equation for a geodesic curve \cite{kolev_lie_2004,cendra_lagrangian_1998}. Indeed in this case, recall that $\omega=\theta|_{\gamma(t)} \dot{\gamma}(t)$ is the left-angular velocity, $\zeta = \omega$ and $\alpha(\omega, \omega)~=~ -\ad^*_\omega(\omega)$. Hence $\gamma$ is a geodesic if and only if $dL_{\gamma(t)}^{-1} \nabla_{\dot \gamma (t)} \dot \gamma(t) = 0$ i.e. setting the left-hand side of \eqref{eq:connection} to $0$. We obtain
\begin{equation}
    \label{eq:ep}
    \begin{cases}
    \dot{\gamma}(t) &= dL_{\gamma(t)} \omega(t) \\
    \dot{\omega}(t) &= \ad_{\omega(t)}^* \omega(t).
    \end{cases}
\end{equation}

\begin{remark}
One can show that the metric is bi-invariant if and only if the adjoint map is skew-symmetric (see \cite{pennec_exponential_2012} or \cite[Prop. 20.7]{gallier_differential_2020}). In this case $ad_\omega^*(\omega) = 0$ and \eqref{eq:ep} coincides with the equation of one-parameter subgroups on $G$.
\end{remark}

\subsection{Reduced Parallel Transport Equation}
The second case is for a vector $Y$ that is parallel along the curve $\gamma$, that is, $\forall t, \nabla_{\dot \gamma(t)} Y(t) = 0$. Similarly to the geodesic equation, we deduce from \eqref{eq:connection} the parallel transport equation expressed in the Lie algebra.

\begin{theorem}
\label{thm}
Let $\gamma$ be a smooth curve on $G$. The vector $Y$ is parallel along $\gamma$ if and only if it is solution to
\begin{equation}
    \label{eq:tp}
    \begin{cases}
    \omega(t) &= dL_{\gamma(t)}^{-1} \dot{\gamma}(t) \\
    Y(t) &= dL_{\gamma(t)} \zeta(t) \\
    \dot{\zeta}(t) &= -\alpha(\omega(t), \zeta(t))
    \end{cases}
\end{equation}
\end{theorem}
Note that in order to parallel transport along a geodesic curve, \eqref{eq:ep} and \eqref{eq:tp} are solved jointly.

\subsection{Application}

We now exemplify Theorem~\ref{thm} on the group  of isometries of $\R^3$, $SE(3)$, endowed with a left-invariant metric $g$. $SE(3)$, is the semi-direct product of the group of three-dimensional rotations $SO(3)$ with $\R^3$, i.e. the group multiplicative law for $R,R' \in SO(3), t,t' \in \R^3$ is given by
\begin{equation*}
    (R,t)\cdot (R',t') = (RR', t + Rt').
\end{equation*}
It can be seen as a subgroup of $GL(4)$ and represented by homogeneous coordinates:
\begin{equation*}
    (R,t) = \begin{pmatrix} R & t \\ 0 & 1 \end{pmatrix},
\end{equation*}
and all group operations then correspond to the matrix operations. Let the metric matrix at the identity be diagonal: $G=\mathrm{diag}(1,1,1,\beta,1,1)$ for some $\beta>0$, the anisotropy parameter. 
An orthonormal basis of the Lie algebra $\se(3)$ is
\begin{align*}
    e_1 = \frac{1}{\sqrt{2}}\begin{pmatrix} 0 & 0 & 0 & 0 \\
            0 & 0 & -1 & 0\\ 0 & 1 & 0 & 0 \\ 0 & 0 & 0 & 0 \end{pmatrix}
    &\qquad&
    e_2 = \frac{1}{\sqrt{2}}\begin{pmatrix} 0 & 0 & 1 & 0 \\
        0 & 0 & 0 & 0\\ -1 & 0 & 0 & 0 \\ 0 & 0 & 0 & 0 \end{pmatrix}
    &\qquad&
    e_3 = \frac{1}{\sqrt{2}}\begin{pmatrix} 0 & -1 & 0 & 0 \\
        1 & 0 & 0 & 0\\ 0 & 0 & 0 & 0 \\ 0 & 0 & 0 & 0 \end{pmatrix}\\
    e_4 = \frac{1}{\sqrt{\beta}}\begin{pmatrix} 0 & 0 & 0 & 1 \\
        0 & 0 & 0 & 0\\ 0 & 0 & 0 & 0 \\ 0 & 0 & 0 & 0 \end{pmatrix}
    &\qquad&
    e_5 = \begin{pmatrix} 0 & 0 & 0 & 0 \\
        0 & 0 & 0 & 1\\ 0 & 0 & 0 & 0 \\ 0 & 0 & 0 & 0 \end{pmatrix}
    &\qquad&
    e_6 = \begin{pmatrix} 0 & 0 & 0 & 0 \\
        0 & 0 & 0 & 0\\ 0 & 0 & 0 & 1 \\ 0 & 0 & 0 & 0 \end{pmatrix} .
\end{align*}
Define the corresponding structure constants $C_{ij}^k = \langle [e_i,e_j],e_k\rangle$, where the Lie bracket $[\cdot,\cdot]$ is the usual matrix commutator. It is straightforward to compute
\begin{align}
    \label{eq:structure_constants}
    C_{ij}^k &= \frac{1}{\sqrt{2}} \;\; \textrm{if} \;\;ijk \;\;\textrm{is a direct cycle of}\;\; \{1,2,3\};\\
    C_{15}^6 &= - C_{16}^5 = - \sqrt{\beta} C_{24}^6 = \frac{1}{\sqrt{\beta}} C_{26}^4 = \sqrt{\beta} C_{34}^5 = -\frac{1}{\sqrt{\beta}} C_{35}^4 = \frac{1}{\sqrt{2}} .
\end{align}
and all others that cannot be deduced by skew-symmetry of the bracket are equal to $0$. The connection can then easily be computed using
\begin{equation*}
    \alpha(e_i, e_j) = \nabla_{e_i}e_j = \frac{1}{2} \sum_k (C_{ij}^k - C_{jk}^i + C_{ki}^j)e_k,
\end{equation*}
For $\beta=1$, $(SE(3), G)$ is a symmetric space and the metric corresponds to the direct product metric of $SO(3) \times \R^3$. However, for $\beta \neq 1$, the geodesics cannot be computed in closed-form and we resort to a numerical scheme to integrate \eqref{eq:ep}.
According to \cite{guigui_numerical_2020}, the pole ladder can be used with only one step of a fourth-order scheme to compute the exponential and logarithm maps at each rung of the ladder.
We use a Runge-Kutta (RK) scheme of order $4$. The Riemannian logarithm is computed with a gradient descent on the initial velocity, where the gradient of the exponential is computed by automatic differentiation. All of these are available in the \url{InvariantMetric} class of the package \href{http://geomstats.ai}{geomstats} \cite{miolane_geomstats_2020}.

\begin{figure}[ht]
    \centering
    \includegraphics[width=12cm]{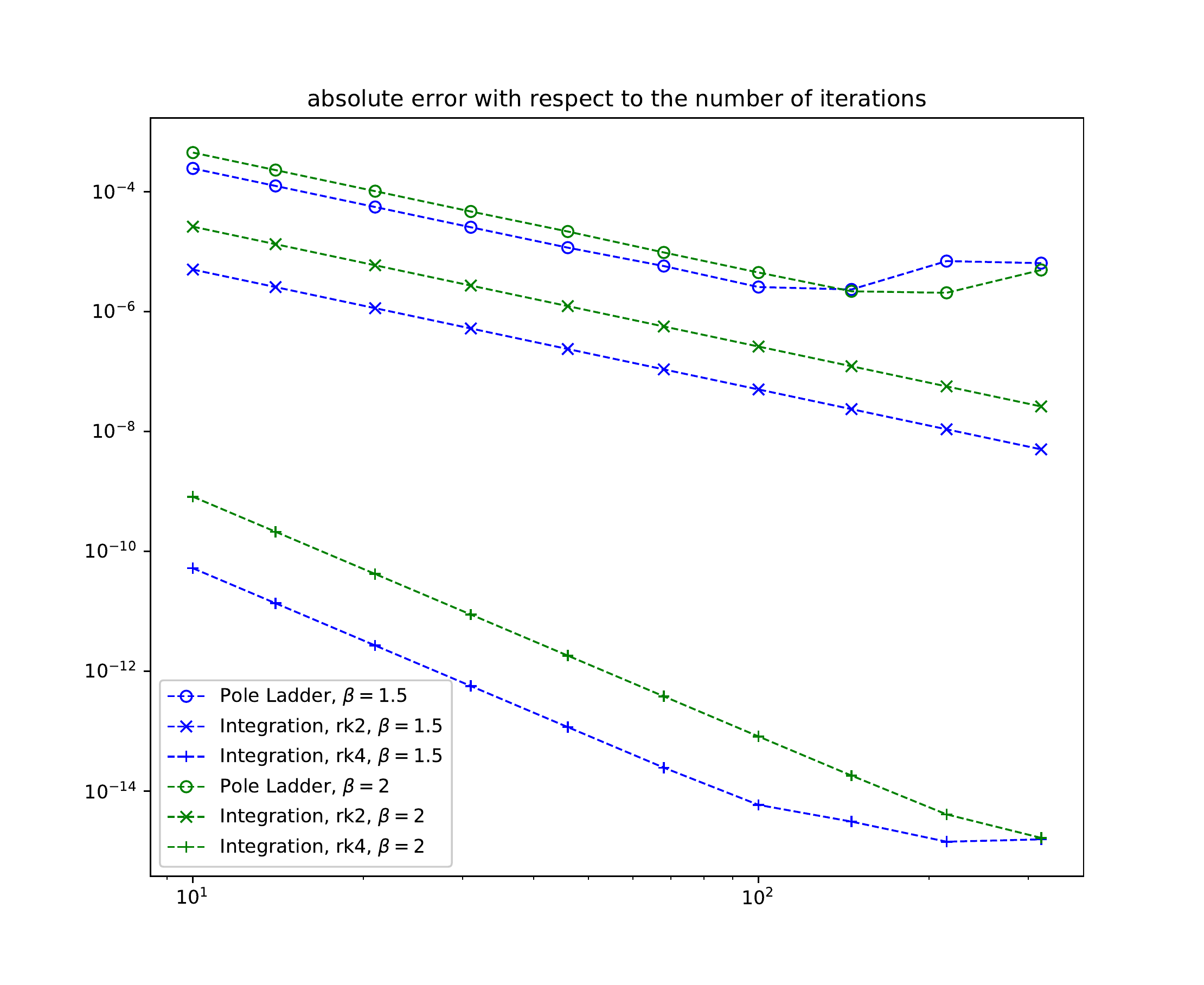}
    \caption{Comparison of the integration of the reduced equation with the pole ladder}
    \label{fig:my_label}
\end{figure}

We now compare the integration of \eqref{eq:tp} to the pole ladder \cite{guigui_numerical_2020} for $\beta=1.5,2$ to parallel transport a tangent vector along a geodesic. The results are displayed on Figure~\ref{fig:my_label} in a log-log plot.
As expected, we reach convergence speeds of order two for the pole ladder and the RK2 scheme, while the RK4 schemes is of order four. Both integration methods are very stable, while the pole ladder is less stable for $\sim n \geq 200$.

\section{Acknowledgments}
\label{sec:acknowledgments}
This work was partially funded by the ERC grant Nr. 786854 G-Statistics from the European Research Council under the European Union’s Horizon 2020 research and innovation program.
It was also supported by the French government through the 3IA Côte d’Azur Investments ANR-19-P3IA-0002 managed by the National Research Agency. 

%
%
%
\bibliographystyle{splncs04}
\bibliography{GSI21_groups.bib}
\end{document}